\documentclass[12pt]{amsart}
\usepackage{hyperref}
\usepackage[capitalise]{cleveref}
\usepackage{xcolor}
\hypersetup{hypertexnames = false,
   bookmarksdepth = 2, bookmarksopen = true, bookmarksnumbered, colorlinks,
  linkcolor={red!50!black},
 citecolor={blue!50!black},
 urlcolor={blue!80!black},
 pdfstartview={XYZ null null 1}}
\usepackage{amsmath,amsthm,epsfig,amssymb}
\usepackage{amsthm,amsfonts,amssymb,amscd}
\usepackage[all]{xy}
\newtheorem{thm}[subsection]{Theorem}

\newtheorem{lem}[subsection]{Lemma}

\theoremstyle{definition}

\newtheorem{proposition-definition}[subsection]{Proposition-Definition}

\newcommand{\CC}{{\mathbb C}}

\newcommand{\ZZ}{{\mathbb Z}}

\newcommand{\PP}{{\mathbb P}}

\newcommand{\OOO}{{\mathcal O}}

\author{Sergey Galkin}
\address{1) PUC-Rio, Departamento de Matemática, Rua Marquês de São Vicente
	225, Gávea, Rio de Janeiro - CEP 22451-900, Brazil.
	2) HSE University, Russian Federation.}
\email{"sergey galkin" <sergey.galkin@phystech.edu> }

\author{D. S. Nagaraj}
\address{Indian Institute of Science Education and Research, 
Rami Reddy Nagar, Karakambadi Road,
Mangalam (P.O.), Tirupati - $517507$,
Andhra Pradesh, INDIA.}
\email{dsn@iisertirupati.ac.in}

\title{ Projective bundles and blow-ups of Projective spaces. }
\date{}
\begin{document}

\begin{abstract}  
The aim of this note is to investigate the relation between two types of 
non-singular projective varieties of Picard rank 2, namely the Projective
bundles over Projective spaces and certain Blow-up of Projective  spaces.
\end{abstract}

\maketitle

\section{Introduction} \setcounter{page}{1} 

Let $X$ be the projective variety obtained from projective space 
$\PP^n= \PP((\CC^{n+1})^*)$ by blowing up along a 
linear subspace $\PP(W^*)$ of co-dimension $r,$  where $\CC^{n+1} \to W$ is
a quotient vector space of dimension $n-r+1.$ If $K$ is the kernel of 
surjective linear map $\CC^{n+1} \to W$ 
resulting projective variety admits a projective bundle 
structure over the smaller dimensional projective space $\PP(K^*).$ 
In fact 
$$ X \simeq \PP(W\otimes\OOO_{\PP(K^*)}\oplus\OOO_{\PP(K^*)}(1)).$$

Question: Are there any other examples of 
blow-ups of a projective space along a sub variety
that also have a structure of 
a projective bundle over a projective space?

Nabanita Ray in \cite{NR} gave examples of non-linear sub varieties in 
$\PP^i, i=3,4,5$  whose blow-ups are
projective bundles over  $\PP^2.$ 

In this note we provide several other examples of  Projective bundles over a Projective space $\PP^n$
which can be realised as a blow up of some $\PP^m$ along a non-linear sub variety.

\begin{thm}\label{thm1}
	Let $n\geq 2$ be an integer and $N= n^2+n-1.$  Let 
	$$E = \underline{Hom}(\OOO_{\PP^n}^n, T_{\PP^n}(-1))$$ 
	be the rank $n^2$ vector bundle on $\PP^n,$ where $T_{\PP^n}(-1)$ is tangent bundle of $\PP^n$ twisted by the inverse of the ample bundle $\OOO_{\PP^n}(1).$
	The projective bundle $\PP(E)$ is isomorphic to a blow-up of 
	$\PP^N=\PP(Hom(\CC^n,\CC^{n+1}))$ 
	along the sub variety of all linear mappings of rank at most $n-1.$
\end{thm}

\begin{thm}\label{thm2} 
	Let $n$ be an even integer and $N= \frac{n(n+1)}{2}-1.$ Let $V= \wedge^2(T_{\PP^n}(-1))$ be the rank $N+1-n= \frac{n(n-1)}{2}$ 
	vector bundle over $\PP^n.$ The projective bundle $\PP(V)$ is isomorphic to  
	$\PP^N=\PP(Alt({n+1}))$ blown-up
	along the sub variety  consists of all alternating linear mappings of rank at most 
	$n-1,$ where $Alt({n+1})$ denote the space of all alternating linear mappings from
	$\CC^{n+1}$ to itself.
\end{thm}

\section{The Projective bundle  $\PP(\underline{Hom}(\OOO_{\PP^n}^n, T_{\PP^n}(-1)))$}

{\bf Proof of Theorem(\ref{thm1}):} On $\PP^n$ by applying
$$\underline{Hom}(\OOO_{\PP^n}^n,-) $$ 
to the standard exact sequence 
\begin{equation}\label{eq01}
 0 \to  \OOO_{\PP^n}(-1) \to 
\OOO_{\PP^n}^{n+1} \to T_{\PP^n}(-1) \to 0.
\end{equation}
we obtain the exact sequence
\begin{equation}\label{eq0}
 0 \to \underline{Hom}(\OOO_{\PP^n}^n, \OOO_{\PP^n}(-1)) \to
\underline{Hom}(\OOO_{\PP^n}^n, \OOO_{\PP^n}^{n+1}) \to E \to 0.
\end{equation}

Note that $E$ is a globally generated vector bundle and ${\rm dim}H^0(E)=n(n+1).$
Hence we obtain a morphism from 
$$\phi : \PP(E) \to \PP^{n(n+1)-1}.$$
 We claim that this is a surjective morphism of degree one. i.e., $\phi$ is a 
 bi-rational morphism and hence $\PP(E)$ is a blow-up of $\PP^{n(n+1)-1}$ 
 [See, Chapter II, Theorem 7.17. cite{Ha}]
  Let $\xi=\OOO_{\PP(E)}(1)$ denote the tautological line bundle on  $\PP(E).$ 
  Then there is an exact
  sequence of vector bundles on $\PP(E)$
\begin{equation}\label{eq1}
 0 \to \Omega^1_{\PP(E)/{\PP^n}} \to \pi^*(E) \to \xi \to 0,
\end{equation}  
where $\pi : \PP(E) \to \PP^n$ is the natural projection.
Since $\xi^{(n(n+1)-1)}$ is the degree of the morphism $\phi$ it is enough to
show that $\xi^{(n(n+1)-1)}=1.$
Since $E$ is vector bundle of rank $n^2$ in the cohomology ring of
$\PP(E)$ the relation 
\begin{equation}\label{eq2}
 \xi^{n^2} = \sum_{i=1}^{n}(-1)^{i+1}\pi^*(C_i(E))\xi^{n^2-i}
\end{equation}
holds, where $C_i(E),\, (1\leq i \leq n)$ are the Chern classes of $E.$
Repeated use of the relation \ref{eq2}  gives 
\begin{equation}\label{eq3}
\xi^{n^2+n-1} = (-1)^n\pi^*(S_n(E))\xi^{n^2-1},
\end{equation} 
where $S_n(E) $ is the $n$th Segre class of $E.$ 
By equation (\ref{eq1}) and the fact that the  total Segre class is the 
inverse of the total Chern class we deduce that 
$S_n(E) = (-1)^n h^n,$
where $h$ is the class of the line bundle $\OOO_{\PP^n}(1)$ in ${\rm H}^2(\PP^n, \ZZ).$
Hence $\xi^{n^2+n-1} = 1.\pi^*(h^n)\xi^{n^2-1}.$ Since $\pi^*(h^n)\xi^{n^2-1}$ is 
the class of a point in $\PP(E)$ it follows that the degree of the map $\phi$ is one.
$\hfill{\Box}$

{\bf Remark:} {\it The bi-rational morphism $\phi$  can be described geometrically
as follows:
   let 
   $$ T = \{ (v, \varphi)\in \PP^n \times \PP({\rm Hom}(\CC^{n},\CC^{n+1}))| (\varphi)^t(v)=0 \},$$
   where $(\varphi)^t :\CC^{n+1} \to \CC^n$ is the transpose of the map
   $\varphi :\CC^{n} \to \CC^{n+1}.$
   Then $T \simeq \PP(E)$  and the bi-rational morphism $\phi$ is the projection
   onto the second factor. Let 
   $$D= \{\varphi \in \PP({\rm Hom}(\CC^{n},\CC^{n+1}))|\,\, \text{rank}(\varphi)\leq n-1\},$$
   $F = \phi^{-1}(D)$ and $U= \PP(E)\setminus F.$
The morphism $\phi$ is one-one on the open set $U.$ }

Next we prove that the set $F$ in the previous Remark is a divisor and we identify
its class in the Picard group of $\PP(E).$ Note that Picard group of $\PP(E)$
is equal to ${\rm H}^2(\PP(E),\ZZ) = \ZZ[\pi^*(h)] \oplus \ZZ [\xi].$

\begin{lem}\label{lem1}
	With the notations of previous Remark the set $F$ is a divisor in $\PP(E)$
	and its divisor class is $ n[\xi]-[\pi^*(h)].$
\end{lem} 

{\bf Proof:}  Rewriting the exact sequence (\ref{eq0}) as 
\begin{equation}\label{eq4}
0 \to \OOO_{\PP^n}(-1)^n \to \OOO_{\PP^n}^{n(n+1)} \to E \to 0
\end{equation}
and taking $n$-th symmetric power
we obtain a long exact sequence
\begin{equation}\label{eq5}
0 \to \OOO_{\PP^n}(-n)\to \OOO_{\PP^n}^{n(n+1)} \otimes \wedge^{n-1}(\OOO_{\PP^n}(-1)^n)\to \cdots 
\end{equation}
\begin{equation}\label{eq6}
\cdots \to S^{n-1}(\OOO_{\PP^n}^{n(n+1)})\otimes \OOO_{\PP^n}(-1)^n \to S^{n}(\OOO_{\PP^n}^{n(n+1)}) \to S^n(E )\to 0.
\end{equation}
Tensoring the long exact sequence (\ref{eq5}) (\ref{eq6}) by $\OOO_{\PP^n}(-1)$ 
 and computing the cohomology
we deduce 
$${\rm H}^0(\PP(E), \xi^{\otimes n}\otimes \pi^*(\OOO_{\PP^n}(-1)))\simeq {\rm H}^0(S^n(E)(-1)) $$
$$ \simeq {\rm H}^n(\OOO_{\PP^n}(-n-1)) \simeq \CC.$$
This proves that the divisor class $ n[\xi]-[\pi^*(h)]$ contains a unique effective divisor. 
From the exact sequence (\ref{eq0}) we deduce that the Segre class $S_{n-1}(E)$ 
is equal to $(-1)^{n-1}nh^{n-1}$ and hence 
$(n[\xi]-[\pi^*(h)]).\xi^{n^2+n-2}= n\pi^*(h^n)\xi^{n^2-1}-(-1)^{n-1}\pi^*(h)\pi^*(S_{n-1}(E))\xi^{n^2-1}= 0.$
For $v \in \PP^n$
the fibre  $E_v$ is isomorphic to ${\rm Hom}(\CC^n, {\CC^{n+1}/{\CC.v}}).$ For
$\varphi \in {\rm Hom}(\CC^n, \CC^{n+1})$  rank of $(\varphi)^t$ is less than or equal to $n-1$ 
if and only if image of $\bar{\varphi}\in {\rm Hom}(\CC^n, {\CC^{n+1}/{\CC.v}})$ of $\varphi$ 
has rank less than or equal to $n-1,$ i.e., $\bar\varphi$ is not an isomorphism. 
 Since the complement of 
isomorphisms in ${\rm Hom}(\CC^n, {\CC^{n+1}/{\CC.v}})$ is given by vanishing of the 
homogeneous polynomial of degree in $n,$ i.e., section of $\OOO_{\PP(E_{v})}(n).$
The set 
$$ F_v =\{\bar{\varphi} \in {\rm Hom}(\CC^n, {\CC^{n+1}/{\CC.v}})|  {\rm rank}(\bar\varphi) \leq n-1 \} $$
is an irreducible divisor given by the vanishing of the restriction to  $\PP(E_v)$ of the 
non-zero section (unique upto non zero scalar) of 
$${\rm H}^0(\PP(E), \xi^{\otimes n}\otimes \pi^*(\OOO_{\PP^n}(-1)))$$
The non-zero section $t$ (unique upto multiplication by a
non zero scalar) determines a section of $\OOO_{\PP(E_{v})}(n).$ It is clear that $F_v$ is the
set of zeros of  the section $t|_{\PP(E_v)}.$  This proves that the set $F$ is the support of the
divisor the $n[\xi]-[\pi^*(h)].$             $\hfill{\Box}$

{\bf Remark 1:} {\it Theorem (\ref{thm1}) can be used to obtain  for 
blow up of $\PP^m, \,(2n-1 \leq m \leq n(n+1)-1),$  along a non linear sub-variety a projective
bundle structure over $\PP^n$ for any integer $n \geq 2.$}

\section{Projective bundle $\PP(\wedge^2(T_{\PP^n}(-1))).$}

{\bf Proof of Theorem(\ref{thm2}):}

  Let  $n \geq 3$ be an integer.  By taking 2nd exterior  power in the exact sequence (\ref{eq01})
  yields the sequence
  \begin{equation}\label{eq7}
 0 \to  T_{\PP^n}(-2) \to 
\wedge^2\OOO_{\PP^n}^{n+1} \to \wedge^2(T_{\PP^n}(-1)) \to 0.
\end{equation}
If we identify ${\rm H}^0(\OOO_{\PP^n}^{n+1})$ with $\CC^{n+1}$ then the 
${\rm H}^0(\wedge^2\OOO_{\PP^n}^{n+1})$ gets identified with $Alt({n+1}),$
where $Alt({n+1})$ denote the set of all alternating linear mappings from
	$\CC^{n+1}$ to itself.  Using the exact sequence (\ref{eq7}) we obtain 
	$$H^0(\wedge^2(T_{\PP^n}(-1))) \simeq Alt({n+1}).$$
	Since $T_{\PP^n}(-1)$ is globally generated $V= \wedge^2(T_{\PP^n}(-1))$  is generated by
$Alt({n+1}).$
 Therefore we get a morphism 
 \begin{equation}\label{eq8}
 \psi : \PP(V) \to \PP(Alt({n+1})).
 \end{equation}
 We claim that $\psi$ is bi-rational i.e., degree of $\psi$ is one and hence
  $\PP(V)$ is a blow-up of $\PP(Alt({n+1})).$ 
 [See, Chapter II, Theorem 7.17. cite{Ha}]
  Let $\zeta=\OOO_{\PP(V)}(1)$ denote the tautological line bundle on  $\PP(V).$ 
  Then there is an exact
  sequence of vector bundles on $\PP(V)$
\begin{equation}\label{eq9}
 0 \to \Omega^1_{\PP(V)/{\PP^n}} \to \pi^*(V) \to \zeta \to 0,
\end{equation}  
where $\pi : \PP(V) \to \PP^n$ is the natural projection.
Since $\zeta^{(\frac{n(n+1)}{2}-1)}$ computes  the degree of the morphism $\psi$ 
 it is enough to
show that, for $n$  even, 
$$\zeta^{(\frac{n(n+1)}{2}-1)}= \pi^*(h^n)\zeta^{\frac{n(n-1)}{2}-1}.$$
 
Since $V$ is vector bundle of rank $n(n-1)/2$ in the cohomology ring of
$\PP(V)$ the relation 
\begin{equation}\label{eq10}
 \zeta^{\frac{n(n-1)}{2}} = \sum_{i=1}^{n}(-1)^{i+1}\pi^*(C_i(V))\zeta^{\frac{n(n+1)}{2}-i}
\end{equation}
holds, where $C_i(V),\, (1\leq i \leq n)$ are the Chern classes of $V.$
Repeated use of the relation \ref{eq10}  gives 
\begin{equation}\label{eq11}
\zeta^{\frac{n(n+1)}{2}-1} = (-1)^n\pi^*(S_n(V))\zeta^{\frac{n(n+1)}{2}-1},
\end{equation} 
where $S_n(V) $ is the $n$th Segre class of $V.$ 
By equation (\ref{eq7}) and the fact that the  total Segre class of $V$ is the 
total Chern class  $C(T_{\PP^n}(-2))$ of $T_{\PP^n}(-2).$ Tensoring the equation (\ref{eq0}) by
$\OOO_{\PP^n}(-1)$ to obtain 
\begin{equation}\label{eq12}
C(T_{\PP^n}(-2))= (1-h)^{n+1}(1-2h)^{-1},
\end{equation}
where $h$ is the class of the line bundle $\OOO_{\PP^n}(1)$ in ${\rm H}^2(\PP^n, \ZZ).$
From the equation (\ref{eq12}) we deduce that 
$$C_n(T_{\PP^n}(-2)) = (\sum_{i=0}^n(-1)^i\binom{n+1}{i}2^{n-i})h^n$$ 
$$=((-1)^{n+1}\frac{(\sum_{i=0}^n\binom{n+1}{i}(-2)^{n+1-i})}{2})h^n$$
$$= ((-1)^{n+1}(\frac{(1-2)^{n+1}}{2})+(-1)^{n+2}\frac{1}{2})h^n.$$
Thus 
$$S_n(V)= \left\{\begin{array}{lcr}
h^n & {\rm if} & n \,\,\, {\rm even}\\
0 & {\rm if} & n \, \,\,{\rm odd}
\end{array} \right.
$$
Now from equation (\ref{eq11}) we deduce that, for $n$ even  
$$\zeta^{\frac{n(n+1)}{2}-1}= \pi^*(h^n)\zeta^{\frac{n(n-1)}{2}-1}.$$
Since $\pi^*(h^n)\zeta^{\frac{n(n-1)}{2}-1}$ is 
the class of a point in $\PP(V)$ it follows that the degree of the map $\psi$ is one
when $n$ is even.
$\hfill{\Box}$
 
 {\bf Remark:} {\it The bi-rational morphism $\psi$  can be described geometrically
as follows:
   let 
   $$ S = \{ (v, \varphi)\in \PP^n \times \PP(Alt({n+1}))| (\varphi)(v)=0 \}.$$
   For $v \in \CC^{n+1}\setminus \{0\}$ the the the linear subspace $\phi \in Alt({n+1})$
   gets identified with fibre  $\wedge^2(T_{\PP^n}(-1))_{v}$ of 
   $\wedge^2(T_{\PP^n}(-1))_{v}$ over the point $[v] \in \PP^n$ hence  
  $S \simeq \PP(V).$   Under this identification, when $n$ is even, the bi-rational morphism 
  $\psi$ is the projection onto the second factor. From now on we assume $n$ is an even
  integer say $n=2k.$ Let 
   $$D= \{\varphi \in \PP(Alt{n+1})|\,\, \text{rank}(\varphi)\leq n-1\},$$
   $F = \phi^{-1}(D)$ and $U= \PP(E)\setminus F.$
The morphism $\psi$ is one-one on the open set $U.$ }

Next we prove that the set $F$ in the previous Remark is a divisor and we identify
its class in the Picard group of $\PP(V).$ Note that Picard group of $\PP(V)$
is equal to ${\rm H}^2(\PP(E),\ZZ) = \ZZ[\pi^*(h)] \oplus \ZZ [\zeta].$

\begin{lem}\label{lem2}
	With the notations of previous Remark the set $F$ is a divisor in $\PP(V)$
	and its divisor class is $ k[\zeta]-[\pi^*(h)],$ where $k =n/2.$   
\end{lem} 

{\bf Proof:}  By equation (\ref{eq7})   
we obtain a long exact sequence for $k$-th symmetric power of $V$
\begin{equation}\label{eq13}
0 \to \wedge^k(T_{\PP^n}(-1))\to \OOO_{\PP^n}^{\frac{n(n+1)}{2}} \otimes \wedge^{k-1}(T_{\PP^n}(-1)))\to \cdots 
\end{equation}
\begin{equation}\label{eq14}
\cdots \to S^{k-1}(\OOO_{\PP^n}^{\frac{n(n+1)}{2}})\otimes (T_{\PP^n}(-1)) \to 
S^{k}(\OOO_{\PP^n}^{\frac{n(n+1)}{2}}) \to S^k(V)\to 0.
\end{equation}
Tensoring the long exact sequence (\ref{eq13}) (\ref{eq14}) by $\OOO_{\PP^n}(-1)$ 
 and computing the cohomology
we deduce 
$${\rm H}^0(\PP(V), \zeta^{\otimes k}\otimes \pi^*(\OOO_{\PP^n}(-1)))\simeq 
{\rm H}^0(S^k(V)(-1)) $$
$$ \simeq {\rm H}^k(\wedge^k(T_{\PP^n}(-1))(-1)) \simeq 
{\rm H}^k((\Omega_{\PP^n}^k)) \simeq \CC,$$
where the $\Omega_{\PP^n}^k$ is the bundle $k$ differential forms and the last isomorphism is the consequence of Bott's formula [See, Page 8 \cite{CMH}].
This proves that the divisor class $ k[\zeta]-[\pi^*(h)]$ contains a unique effective divisor. 
From the exact sequence (\ref{eq7}) we deduce that the Segre class $S_{n-1}(V)$ 
is equal to $(-1)^{n-1}kh^{n-1}$ and hence 
$$(k[\zeta]-[\pi^*(h)]).\zeta^{\frac{n(n+1)}{2}-2}= $$
$$\pi^*(h^n)\zeta^{\frac{n(n-1)}{2}-1}-(-1)^{n-1}\pi^*(h)\pi^*(S_{n-1}(V))\zeta^{\frac{n(n+1)}{2}-1}= 0.$$
For $v \in \PP^n$
the fibre  $V_v$ is isomorphic to ${\rm Alt}(\CC^n, {\CC^{n+1}/{\CC.v}}).$ For
$\varphi \in Alt{n+1}$  rank of $\varphi$ is less than or equal to $n-1$ 
if and only if image of $\bar{\varphi}\in {\rm Alt}(\CC^n, {\CC^{n+1}/{\CC.v}})$ of $\varphi$ 
has rank less than or equal to $n-1,$ i.e., $\bar\varphi$ is not an isomorphism. 
 Since the complement of 
isomorphisms in ${\rm Alt}(\CC^n, {\CC^{n+1}/{\CC.v}})$ is given by vanishing of the 
homogeneous polynomial of degree in $k,$ i.e., section of $\OOO_{\PP(V_{v})}(k).$ 
The set 
$$ F_v =\{\bar{\varphi} \in {\rm Alt}(\CC^n, {\CC^{n+1}/{\CC.v}})|  {\rm rank}(\bar\varphi) \leq n-1 \} $$
is an irreducible divisor given by the vanishing of the restriction to  $\PP(V_v)$ of the 
non-zero section (unique upto non zero scalar) of 
$${\rm H}^0(\PP(V), \zeta^{\otimes k}\otimes \pi^*(\OOO_{\PP^n}(-1)))$$
The non-zero section $t$ (unique upto multiplication by a
non zero scalar) determines a section of $\OOO_{\PP(V_{v})}(k),$ namely  Pfaffian
whose square is the determinant of skew symmetric form. It is clear that $F_v$ is the
set of zeros of  the section $t|_{\PP(V_v)}.$  This proves that the set $F$ is the support of the
divisor the $k[\zeta]-[\pi^*(h)].$             $\hfill{\Box}$

{\bf Remark 2:} {\it Theorem (\ref{thm2}) can be used to obtain  for 
blow up of $\PP^m, \,(2n-1 \leq m \leq \frac{n(n-1)}{2}-1)$   along a non linear sub-variety a projective
bundle structure over $\PP^n,$ for even integer $n \geq 4.$}

{\bf Remark 3:} The special case of Theorem (\ref{thm2}) appears first in \cite{EL} and has
been used in the context of Quantum Cohomology in \cite{CCGK}, \cite{CGKS}, \cite{AS1},
\cite{AS2}, \cite{AS3}.

{\bf Acknowledgement: } We would  like to thank organisers of the  conference 
on "Derived Category and Algebraic Geometry" at TIFR Mumbai for the invitation to give a talk and TIFR for its hospitality.

\providecommand{\arxiv}[1]{\href{http://arxiv.org/abs/#1}{arXiv:#1}}
\providecommand{\doi}[1]{\href{http://doi.org/#1}{\tt doi:#1}}

\end{document}